\documentclass[hidelinks,onefignum,onetabnum]{siamart220329}

\usepackage{lipsum}
\usepackage{amsfonts}
\usepackage{graphicx}
\usepackage{epstopdf}
\usepackage{algorithmic}
\ifpdf
  \DeclareGraphicsExtensions{.eps,.pdf,.png,.jpg}
\else
  \DeclareGraphicsExtensions{.eps}
\fi

\newsiamremark{remark}{Remark}
\newsiamremark{hypothesis}{Hypothesis}
\crefname{hypothesis}{Hypothesis}{Hypotheses}
\newsiamthm{claim}{Claim}

\headers{An observer 
for 
pipeline flow with hydrogen blending}{M. Gugat and J. Giesselmann}

\title{An observer for 
pipeline flow with hydrogen blending
in gas networks: Exponential synchronization
\thanks{Submitted to the editors 04/05/2023.
\funding{
This work was funded by DFG in the 
 Collaborative Research Centre
CRC/Transregio 154,
Mathematical Modelling, Simulation and Optimization Using the Example of Gas Networks,
Project C03 and   C05.
The authors thank the Bundesministerium für Bildung und Forschung (BMBF)
for funding under DAAD grant 57654073 'Uncertain data in control of PDE systems'.
}}}

\author{Martin Gugat\thanks{
Friedrich-Alexander-Universit\"at Erlangen-N\"urnberg (FAU),
Department Mathematik, Cauerstr. 11, 91058 Erlangen, Germany, 
(\email{martin.gugat@fau.de})
}
\and
Jan Giesselmann
\thanks{Technische Universit\"at Darmstadt,
Fachbereich Mathematik, Dolivostr. 15, 64293 Darmstadt, Germany.
(\email{giesselmann@mathematik.tu-darmstadt.de})
}
}

\usepackage{amsopn}

\ifpdf
\hypersetup{
  pdftitle={An observer for pipeline flow with hydrogen blending},
  pdfauthor={M. Gugat,  and J. Giesselmann}
}
\fi

\usepackage{float}

\usepackage{yfonts}

\usepackage{amsfonts,amsmath,amssymb}
\usepackage{graphicx}

\renewcommand{\geq}{\geqslant}
\renewcommand{\leq}{\leqslant}

\usepackage{siunitx}

\usepackage{color}

\def\martin#1{{\color{black}#1}}

\newcommand{\tMz}{\mbox{$\hat M$} }

\begin{document}

\maketitle




\begin{abstract}
 \textcolor{black}{
We consider a state estimation problem for gas flows in pipeline networks where 
hydrogen is blended into the natural gas. 
The flow is modeled by the quasi-linear 
isothermal 
Euler equations coupled to an advection equation on a graph. 
The flow through the  vertices  where the pipes are connected
is governed by algebraic node conditions.
The state is approximated by an observer system that uses nodal measurements.
We prove that the
state of the observer system
converges to the original system state
exponentially fast in the $L^2$-norm if the measurements are exact. 
If measurement errors are present we show that the observer state approximates the original system state up to an error that is proportional to the maximal measurement error.
The proof of the synchronization result uses
Lyapunov functions with exponential weights.
}
\end{abstract}









\section*{Introduction}
Hydrogen is a promising condidate as an energy carrier in a decarbonized economy. 
In order to use existing infrastructure for natural gas transport 
hydrogen blending is of interest.
We present an observer system for the flow of gas through
networks of pipelines that is governed by the 
isothermal
Euler equations
with an additional conservation law to allow for hydrogen blending.
We consider classical solutions
that allow to investigate the isothermal Euler equations
in  quasilinear form.
At the nodes of the  gas pipeline network the solutions
for the adjacent pipes are coupled by  algebraic node conditions that require
the conservation of mass and the continuity of the pressure.
Moreover, we assume perfect mixing of the blended hydrogen
and natural gas. 

We present a nodal  observer
system for the blended flow through the network
where the coupling to
the original system 
is governed 
at each 
vertex 
$v$ of the graph 
by a real  parameter $\mu^v\in [0, 1]$.
We show that if 
$\mu^v\geq 0$ is sufficiently small,
the observer system yields an approximation of the
full state in the original system 
where
in the absence of measurement 
{errors}
the $L^2$-norm  of the observer error decays  exponentially fast
with time 
if the $C^1$-norm of the initial state is sufficiently small
and generates a semi-global classical solution.
This assumption makes sense since in 
the operation of gas pipelines,
the velocity of the gas flow is usually quite small
compared with the speed of sound.
If there is a pointwise upper bound
for the measurement error,
the error in the state estimation obtained
from the observer decays exponentially fast
up to the level of the measurement error.



Hence 
the proposed observer yields 
a reliable approximation for the complete
state in the original system exponentially fast. 
%
%
Similarly as in \cite{cobada} 
the proof of the exponential synchronization 
uses Lyapunov functions with exponential weights.

%
%
For semi-linear hyperbolic equations
observers based on  distributed measurements have been presented in
\cite{CCDB12,CIM15} and for quasi-linear hyperbolic equations in \cite{BMP15,GiesselmannGugatKunkel}
A boundary observer for semi-linear hyperbolic problems
based on  the backstepping method is  studied \cite{HAK16}.
In \cite{gugiku} an observer for semi-linear hyperbolic equations based upon
nodal measurements is presented.
The proofs in \cite{gugiku} are based upon a suitable observability  inequality. 
The recovery of an unknown initial state
using an observer 
is studied in \cite{IMT20}.
The design of boundary observers  for hyperbolic PDEs
coupled with 
ODEs 
is studied
in \cite{FERRANTE}.

The novelty of our contribution is the construction of 
an observer for a  system that is governed by 
 networked quasilinear PDEs and uses 
observations that are given on isolated points in  space only.
In addition, we include the effect of measurement errors in
our analysis.
Motivated by current discussions in decarbonization,
our model includes the case of blended flow with
a small concentration of hydrogen.





%

This paper has the following structure.
In Section \ref{isothermalEuler} we introduce
the 
isothermal Euler equations
with an additional conservation
law as a model for hydrogen blending.
In Section  \ref{sec:Riemann} we present the
corresponding Riemann invariants and  transform the system
in diagonal form.
In Section \ref{sec:networks} we present
the node conditions that model the flow 
of the gas mixture through
the junctions in  a gas  pipeline network.
In Section \ref{systems},
we state the  observed system  \textbf{(S)}, the observer system \textbf{(R)}
and  the system \textbf{(Diff)}  that is satisfied by the observation error.
 For our analysis, we use the framework of semi-global classical solutions
that are defined through a fixed point iteration
along the characteristic curves.
A well--posedness result 
is presented in 
Section \ref{wellposedness}.
It is of a similar type as the results 
in \cite{libook}.

In Section \ref{l2lyapunov} the 
decay 
of the $L^2$-error of the observer 
for the density and the flow rate of the mixture is studied 
under smallness assumptions for the initial state.
We show that if there is no measurement error it 
converges  to zero
exponentially fast.
In the proof a quadratic $L^2$-Lyapunov function
with exponential weights is used.
Also in this analysis, possible measurement error is included.
It is shown in the presence of measurement error
the observer error
converges exponentially fast to
a level that is of the order of the measurement error.
The proof of this statement is based upon an 
appropriate version of Gronwall's Lemma.
In Section \ref{l2lyapunovh} 
we show the exponential decay of
the observer error for hydrogen concentration.
The result is of a similar type as in 
Section \ref{l2lyapunov}. Note that 
in  the proof of the result in
Section \ref{l2lyapunovh}  
the result from Section \ref{l2lyapunov} is 
also used.




\section{The isothermal Euler equations}
\label{isothermalEuler}
\label{motivation}
\noindent
The one-dimensional isothermal Euler equations are a well-known 
 model for the flow of gas  through pipelines,
 see for example   \cite{bandahertyklar}.
A pipeline network corresponds to 
a  finite  graph $G=(V,\,E)$ where 
$V$ denotes the set of vertices and $E$ the set of edges.
Each pipe of length $L^e>0$ is parameterized by the interval $[0,\, L^e]$ 
and  corresponds to an edge $e\in E$.
The diameter of the pipe is denoted by 
$D^e>0$  and
the  friction coefficient by 
$\lambda_{fric}^e>0$.
Let 
$\theta^e := \frac{\lambda_{fric}^e}{D^e}$.
%
For the density
of the gas mixture that flows through  pipe $e$ 
we use the notation 
$\hat \rho^e$,
for  the pressure  $\hat p^e$
 and for  the mass flow rate $\hat q^e$.
The isothermal Euler equations are hyperbolic provided 
that the pressure $\hat p= p(\hat \rho)$ is given as a monotone increasing function of the density. We  assume that this monotonicity is strict.
Examples are the isentropic gas law $p(\hat \rho)= a(\hat \rho)^\gamma$ with $a >0$, $\gamma > 1$ and the model of the American Gas Association (AGA), see \cite{gugatulbrich}, 
$p(\hat \rho)= \tfrac{ R_s \,{\cal T} \hat \rho}{1- \tilde \alpha \hat \rho}$ 
where ${\cal T}$ denotes the temperature, 
$R_s$  is the gas constant and $\tilde \alpha \leq 0$.
For $\tilde \alpha=0$ the AGA model reduces to the ideal gas  law $p(\hat \rho)=R_s \, {\cal T} \hat \rho$.
We study a model that is based upon the $2 \times 2$ Euler equations
\begin{equation}
\label{isothermaleuler}
\left\{
\begin{array}{l}
{\hat \rho}_t^e + {\hat q}_x^e = 0,
\\
{\hat q}_t^e + \left( {\hat p}^e + \frac{({\hat q}^e)^2}{{\hat \rho}^e} \right)_x = 
- \frac{1}{2} \theta^e \frac{{\hat q}^e \; |{\hat q}^e|}{{\hat \rho}^e}
\end{array}
\right.
\end{equation}
that govern the flow of the mixture through  pipe $e$.
Similarly as in \cite{zlotnik}, let 
$\rho_{(h)}^e$ denote the density
of the hydrogen that is blended in the mixture 
and
$q_{(h)}^e$ the corresponding  mass flow rate.
We assume that
$  q_{(h)}^e = W(\hat q^e,\, \hat \rho^e) \,\rho_{(h)}^e $.
This means that the velocity 
$W(\hat q^e,\, \hat \rho^e)$
of the hydrogen is
determined by the flow of the gas mixture.
This assumption is justified if the concentration of
the hydrogen is sufficiently low.
We assume that 
$W = \frac{\hat q^e}{\hat \rho^e + \gamma}$ 
where $\gamma \in  [0,\, \infty)$.
For $\gamma =0$
this means that the velocity of the hydrogen is
equal to the velocity of the gas mixture.
A choice of $\gamma>0$ allows to model
a situation where the hydrogen velocity
is smaller than the velocity of the mixture. 
For all $e\in E$, the 
conservation of mass yields the additional 
conservation law 
 \begin{equation}
     \label{hydrogen}
  (\rho_{(h)}^e)_t
 + 
(q_{(h)}^e
)_x=0.
 \end{equation}\,


\subsection{The system in diagonal form}
\label{sec:Riemann}
As stated in \cite[Chapter 7.3]{dafermosbook} for  every $2 \times 2$ system of hyperbolic conservation laws
(for example for  \eqref{isothermaleuler})
we can find a system of Riemann invariants.
Although the extended 
 system \eqref{isothermaleuler}, \eqref{hydrogen}  is a $3\times 3$ system, it is endowed with
 Riemann invariants that allow to write the
system in diagonal form with the
Riemann invariants as new variables.
For this purpose we introduce the notation  
$\;\tilde R(\rho)=    \int_1^\rho \frac{\sqrt{ p'}(r) }{r}   \, dr
$.
We obtain the Riemann invariants
$ R_\pm (\hat \rho, \hat q)
=  \tilde R(\hat \rho) \pm
\frac{\hat q}{\hat \rho} 
$
and
$R_0 
= \frac{ \rho_{(h)} } {  \hat \rho +  \gamma }.
$
With the eigenvalues 
$\lambda_\pm = \frac{\hat q}{\hat \rho}  \pm \sqrt{ p' }$
and 
$\lambda_0 = W(\hat q, \, \hat \rho)$  
we can write the system
\eqref{isothermaleuler}, \eqref{hydrogen}
in the diagonal form 
\begin{multline}
    \label{hydrodiag}
\begin{pmatrix}
R_+^e
\\
R_-^e
\\
R_0^e 
\end{pmatrix}_t
+
\begin{pmatrix}
\lambda_+(R_+^e, R_-^e) & 0 & 0
\\
 0 & \lambda_-(R_+^e, R_-^e) &
0
\\
 0 &   0 & \lambda_0(R_+^e, \, R_-^e)
\end{pmatrix}
\begin{pmatrix}
R_+^e
\\
R_-^e
\\
R_0^e 
\end{pmatrix}_x
\\
=
\begin{pmatrix}
- \tfrac{\theta^e}{8} |R_+^e - R_-^e| \, (R_+^e - R_-^e)
\\
 \tfrac{\theta^e}{8}   |R_+^e - R_-^e| \, (R_+^e - R_-^e)
\\
0
\end{pmatrix}
.
\end{multline}
Here the eigenvalues are represented as functions of
the Riemann invariants, e.g. 
\[
 \lambda^e_\pm
=
\frac{R_+^e - R_-^e}{2} 
\pm 
\sqrt{p'\left(\tilde R^{-1}   \left(\frac{R_+^e + R_-^e}{2}\right) \right)  }
.
\]
The case for the AGA model for real gas  is presented in detail in  \cite{gugatulbrich}.

\section{The Node Conditions}
\label{sec:networks}
 In our model,
 the flow through the nodes of the network
 is governed by algebraic coupling conditions.
%
%
For any node $v \in V$ 
let $E_0(v)$ denote
the set of edges in the graph that are incident to $v$.
Let $x^e(v)\in \{0,L^e\}$ denote the
end of the interval~$[0,L^e]$ that corresponds
to the edge $e \in E_0(v)$.
Define
\begin{equation}
\label{mathfraksdefinition}
{\textbf n}(v, \, e):=\left\{
\begin{array}{rll}
-1 & {\rm if} & x^e(v)=0 \;\mbox{\rm and }\; e\in E_0(v),
\\
1 & {\rm if} & x^e(v)=L^e \;\mbox{\rm and }\; e\in E_0(v),
\\
0 &{\rm if}  &  e\not\in E_0(v).
\end{array}
\right.
\end{equation}
The conservation of mass at the interior nodes
(i.e. $v\in V$ with  $|E_0(v)|\geq 2$)  
is represented by  the 
Kirchhoff condition
\begin{equation}
\label{Kirchhoff}
\sum_{e \in E_0(v)}
{\textbf n}(v, \, e)
\,
 (D^e)^2\,
\hat q^e(x^e(v))
=  
0
.
\end{equation}
As a further coupling condition, we  require the  continuity of the pressure at $v$. This choice leads to well-posed Riemann problems, see \cite{bandahertyklar}. 
It  means that
 \begin{equation}
\label{pressurecontinuity}
p(\hat \rho^e(t,x^e(v))) = p(\hat \rho^f(t,x^f(v)))
\; \mbox{
for all $e$, $f\in E_0(v)$.
}
\end{equation}
Another possible choice is advocated by \cite{reigstad}, namely the continuity of enthalpy:
 \begin{equation}
\label{enthalpycontinuity}
F'(\rho^e(t,x^e(v)))  + \frac{(q^e(t,x^e(v))^2}{(\rho^e(t,x^e(v))^2} = F'(\rho^e(t,x^f(v)))
+ \frac{(q^f(t,x^e(v))^2}{(\rho^f(t,x^e(v))^2}
\end{equation}
\mbox{
for all $e$, $f\in E_0(v)$
}
where $F(\rho)$ is the pressure potential that is defined  by
$
F(\rho) = \rho \int_1^\rho \frac{p(r)}{r^2} dr.
$
 In the small velocity limit 
 $\tfrac{q}{\rho} \rightarrow 0$  both \eqref{pressurecontinuity} and \eqref{enthalpycontinuity} enforce the continuity of densities, since $p(\rho)$ and $F'(\rho)$ are both strictly monotone increasing. 
 Thus both conditions coincide asymptotically in the low velocity limit.


Now we state the node conditions 
\eqref{Kirchhoff}, \eqref{pressurecontinuity}  
in terms of the Riemann invariants $(R_+,\, R_-)$.
Similarly  as in  \cite{gugiku}, 
define the components of the vectors 
$R_{in }^v(t)$,
$R_{out}^v(t)
\in {\mathbb R}^{|E_0(v)|}
$
as follows: For $e  \in  E_0(v)$ we set 
\begin{equation}
\label{indefinition}
R_{in }^e(t,  \, x^e(v))
:=\left\{
\begin{array}{rll}
R^e_-(t,\, x^e(v))  & {\rm if} & x^e(v)=0 \;\mbox{\rm and }\; e\in E_0(v),
\\
R^e_+(t,\, x^e(v)) & {\rm if} & x^e(v)=L^e \;\mbox{\rm and }\; e\in E_0(v).
\end{array}
\right.
\end{equation}
\begin{equation}
\label{outdefinition}
R_{out}^e(t, \, x^e(v))
:=\left\{
\begin{array}{rll}
R^e_+(t,\, x^e(v))  & {\rm if} & x^e(v)=0 \;\mbox{\rm and }\; e\in E_0(v),
\\
R^e_-(t,\, x^e(v)) & {\rm if} & x^e(v)=L^e \;\mbox{\rm and }\; e\in E_0(v).
\end{array}
\right.
\end{equation}

For $v \in V$ and $e\in E_0(v)$ 
with 
$|E_0(v)| \geq 2$
define
$ \omega_v :=  \tfrac{2}{\sum_{ f \in E_0(v) } (D^f)^2 }$
and
\[K^{v,e}_\sigma(R_+, R_-,t)=
 - 
 R^e_{in}(t,  \, x^e(v))  
 +
 \omega_v \sum\limits_{g\in E_0(v)} (D^g)^2 \, 
    R^g_{in}(t,  \, x^g(v)).
    \]

For the convenience of the reader we state 
the following lemma 
from \cite{gugiku}:
\begin{lemma}
\label{lemmaeins}
For any node $v\in V$ with 
$|E_0(v)| \geq 2$ 
and \mbox{ $e\in E_0(v)$}
 the node conditions
(\ref{Kirchhoff}),
(\ref{pressurecontinuity})
can be written in terms of the vectors 
$R_{in }^v(t)$,
$R_{out}^v(t)$
as 
\begin{equation}\label{Omegav}
R_{out}^e(t,  \, x^e(v)) = 
K^{v,e}_\sigma(R_+, R_-,t).
\end{equation}
\end{lemma}
The proof of Lemma \ref{lemmaeins} is given in  \cite{gugiku}.

To complete our model,
we also impose a coupling condition 
for the hydrogen flow.
The conservation of mass requires the
Kirchhoff condition
\begin{equation}
\label{KirchhoffHydrogen}
\sum_{e \in E_0(v)}
{\textbf n}(v, \, e)
\,
 (D^e)^2\,
q^e_{(h)}(t,\, x^e(v))
=  
0
.
\end{equation}
For $v\in V$ and $t\geq 0$  define the set
\[
E_{in}(v,t)
=
\{e \in E_0(v):\,
\textbf{n}(v,\, e) \,\lambda_0(R^e_+(t,\, x_e(v)),\, R^e_-(t,\, x_e(v))) \geq 0
\}
\]

and  require
for 
all
$e\in E_0(v)$
with
$e\not \in E_{in}(v,t)$
the perfect mixing condition 
(see \cite{colombomauri},
\cite {colombogaravello}) 
\begin{equation}
\rho_{(h)}^e(t,\, x^e(v))
=
\sum_{f\in E_{in}(v,t)}
 \lambda^f_R(t) \, \rho_{(h)}^f(t,\, x^f(v))
\end{equation}
where 
\begin{equation}
\label{lambdafdefinition}
\lambda^f_R(t)  = \tfrac{(D^f)^2 |\lambda_0^f| (t, x^f(v)) }{\sum\limits_{g\in E_{in}(v,t) }    (D^g)^2 |\lambda_0^g|(t, x^g(v))   }
\end{equation}
and, by misuse of notation, we write $\lambda_0(t,x)$ instead of $\lambda_0(R_+(t,x),R_-(t,x))$.
Note that perfect mixing implies 
(\ref{KirchhoffHydrogen}) 
For $v \in V$ with 
$|E_0(v)| \geq 2$
let
\begin{equation}
    \label{k0definition}
K^{v,e}_0(R_0, \, R_+, \, R_-,\, t)=
\sum\limits_{f\in E_{in}(v,\,t)}
\lambda^f_R(t) 
\,  R_0^f(t,\, x^f(v))
\end{equation}
if $e\not \in E_{in}(v,t)$.
%
For the Riemann invariant $R_0$ this yields
\begin{equation}
    \label{R0Bedingung}
R_0^e(t,\, x^e(v))
=
K^{v,e}_0(R_0, \, R_+, \, R_-,\, t)
 \mbox{ if } e\not \in E_{in}(v,t).
 \end{equation}

For a boundary node $v\in V$ where $|E_0(v)|=1$ 
we state the boundary  conditions
for $R_+^e$, $R_-^e$ respectively
in the form
\begin{eqnarray}
\label{rbriemann2neu}
R^e_{out}(t,\, x^e(v)) & = 
( 1 - \mu^v) \,  u_{\sigma}^v(t) + \mu^v R^e_{in}(t,\, x^e(v))
\end{eqnarray}
where $\mu^v\in 
[0, 1]$
is a given number.
If 
$e\not \in E_{in}(v,t)$
also a boundary condition for $R_0$ is required.
We state it as
\begin{eqnarray}
\label{rbriemann0}
R^e_0(t,\, x^e(v)) & = 
 u^e_0(t).
\end{eqnarray}
The regularity requirements for
the functions
$ u^e_{\varsigma} $ and $u^e_0 $
will be stated in the section about the 
well-posedness of the system.

\begin{remark}
For $\mu^v=0$, \eqref{rbriemann2neu} is a Dirichlet boundary condition for the incoming Riemann invariant. 
For $\mu^v=1$, 
it 
is a Dirichlet boundary condition for the velocity.


\end{remark}


\section{The observed system \textbf{(S)}, the observer system \textbf{(R)}
and the error system \textbf{(diff)}}
\label{systems}
In this section we state the  observed system,
 we introduce the observer system and 
 state the system that is satisfied by the observation error.


For $e\in E$, let $\nu^e   =  \frac{1}{8}\theta^e $ 
be given and define
\begin{equation}
\label{sigmadefinition}
\sigma^e (R_+^e,\,R_-^e)
=
\nu^e  \, \left|R_+^e - R_-^e\right| \, (R_+^e - R_-^e).
\end{equation}
Define $ \hat \Delta^e(R_+, R_-)$ as the  diagonal $2 \times 2$ matrix
that contains the eigenvalues $\lambda_+$ and $\lambda_-$
 corresponding to the edge $e\in E$.
In terms of the Riemann invariants, the quasilinear  system (\ref{isothermaleuler}) has the  diagonal form 
\begin{equation}
    \label{pqdiag}
\partial_t 
\left(
\begin{array}{r}
R_+^e
\\
R_-^e
\end{array}
\right)
 +
\hat  \Delta^e(R_+, R_-) \,
\partial_x
\left(
\begin{array}{r}
R_+^e
\\
R_-^e
\end{array}
\right)
 =
\sigma^e (R_+^e,\,R_-^e)
\;
\left(
\begin{array}{r}
-1
\\
1 \end{array}
\right).
\end{equation}
This system is independent of $R_0$.
Hence 
the $3\times 3$ system (\ref{hydrodiag}) can
be solved by solving
the $2\times 2$ system 
(\ref{pqdiag}) first and then in a second step 
the equation for $R_0$ with the
given values for $(R_+, R_-)$.
This approach is possible since also
the node conditions for
$(R_+. R_-)$ are independent of $R_0$.



%
%
Define $ \Delta^e(R_+, R_-)$ as the  diagonal $3 \times 3$ matrix
that contains the eigenvalues $\lambda_+$, $\lambda_-$
and $\lambda_0$ corresponding to the edge $e\in E$.
The  quasilinear model with hydrogen blending has the following form:
\[
{\bf (S)}
\left\{
\begin{array}{l}
S_\iota^e(0, x)     = y_\iota^e (x),
\iota \in \{+, \, -,  \, 0 \}, \,  x\in (0,L^e),  \, e\in E,
\\
\\
S_{out}^e(t, x^e(v))   = 
(1- \mu^v) \, u_{\sigma}^v(t)   + \mu^v S_{in}^e(t,  x^e(v) )  , \,  t\in (0,T),\;  
\\ 
\mbox{\rm if } 
 e\in E_0(v) \mbox{  and  } \;|E_0(v)|=1;
\\ 
S_0^e(t, x^e(v))   =   u^e_0(t), \,  t\in (0,T),\;  
\\
 \mbox{\rm if } 
 e\in E_0(v), \; |E_0(v)|=1
 \mbox{  and  }
e\not \in  E_{in}(v,t);
\\
\\
S_{out}^e(t, x^e(v))     = K^{v,e}_\sigma(S_+, S_-,t)
 ,\,  t\in (0,T), 
 \\
    \mbox{\rm if } \; e\in E_0(v) \mbox{  and  } \;|E_0(v)| \geq 2;
\\
S_{0}^e(t,\, x^e(v))
=
K^{v,e}_0(S_0, \, S_+, \, S_-,\, t)
 \\ 
  \mbox{\rm if }
  e\in E_0(v), \,  |E_0(v)| \geq 2 \,  \mbox{  and  }\, e \not \in E_{in}(v, t);
\\
\\
 \partial_t\!
\left(
\begin{array}{r}
S_+^e
\\
S_-^e
\\
S_0^e
\end{array}
\right)
 +
  \Delta^e(S_+, S_-)\,
\partial_x\!
\left(
\begin{array}{r}
S_+^e
\\
S_-^e
\\S_0^e
\end{array}
\right)
 =
\sigma^e (S_+^e,\,S_-^e)
\;
\left(
\begin{array}{r}
-1
\\
1
\\
0
\end{array}
\right)
\;\;
\\
\mbox{\rm on }\; [0,T]\times [0,L^e],\, e\in E.
\end{array}
\right.
\]



Now we introduce the  observer system  ${\bf (R)} $
that depends on numbers
$\mu^v\in  [0,1]$
that are given for all $v\in V$
and control
the flow of information from
the original system to the observer system.
For an interior node with $\mu^v=0$,
the values 
at the node $v$ in the observer system
are fully determined by
the information from the  system ${\bf (S)} $.
For $\mu^v=1$, 
no data from ${\bf (S)} $
are used at the node.
We study the case that
for the data from ${\bf (S)} $ 
we have  perturbations
that is smoothed to
give a $C^1$-function of time
$Z=(Z_+, \, Z_-,\, Z_0)$
for all $v\in V$, $e\in E_0(v)$
at the points $x^e(v)$.
%
The observer ${\bf (R)}$ is defined as follows:
\[
{\bf (R)}
\left\{
\begin{array}{l}
R_\iota^e(0, x)     =  z_\iota^e (x),
\iota \in \{+, \, -,  \, 0 \}, \,  x\in (0,L^e),  \, e\in E,
\\
\\
R_{out}^e(t, x^e(v))   =  (1-\mu^v) [ u^v_{\sigma}( t )  + Z_{out}^e(t, x^e(v)) ]
+ \mu^v  \, R_{in}^e(t, x^e(v)),
\\
  t\in (0,T),
  \mbox{\rm if } \;
  e\in E_0(v) \mbox{  and  }
  |E_0(v)|=1;
\\ 
 %
 R_0^e(t, x^e(v))   =   [u^e_0(t) + Z_{0}^e(t, x^e(v)) ], 
 \\
 t\in (0,T),\,
 \mbox{\rm if } 
 e\in E_0(v), \; |E_0(v)|=1
 \mbox{  and  }
e\not \in  E_{in}(v,t);
\\
\\
  R^e_{out}(t,\, x^e(v))=
  \\
  \hfill 
  \mu^v\, K^{v,e}_\sigma(R_+, R_-,t) 
  +
  (1 - \mu^v) \,[ S^e_{out}(t,\, x^e(v)) + Z_{out}^e(t, x^e(v)) ],
\\
t\in (0,T), 
 \,  \mbox{\rm if } \;
 e\in E_0(v) \mbox{  and  } 
  |E_0(v)|\geq 2;
\\
R_{0}^e(t,\, x^e(v))
=
 \\
  \hfill 
 \mu^v\, K^{v,e}_0(R_0, \, R_+, \, R_-,\, t)
 +
   (1 - \mu^v) \,[ S^e_{0}(t,\, x^e(v))+ Z_{0}^e(t, x^e(v)) ],
 \\ 
t\in (0,T),  \, \mbox{\rm if }
  e\in E_0(v), \,  |E_0(v)| \geq 2 \,  \mbox{  and  }\, e \not \in E_{in}(v, t);
\\
\\
\partial_t\!
\begin{pmatrix}
R_+^e
\\
R_-^e
\\
R_0^e
\end{pmatrix}
 +
  \Delta^e(R_+, R_-)\,
\partial_x\!
\begin{pmatrix}
R_+^e
\\
R_-^e
\\
R_0^e
\end{pmatrix}
 =
\sigma^e (R_+^e,\,R_-^e)
\,
\begin{pmatrix}
-1
\\
1
\\
0
\end{pmatrix}
\\
\mbox{\rm on }\, [0,T]\times [0,L^e],\, e\in E.
\end{array}
\right.
\]
The initial state
$(z_+^e, z_-^e, z^e_0)$
is an approximation  of the initial state of the original system.
The data  from the original system ${\bf (S)}$ 
enters the system state through the node conditions.
Indeed, for $\mu^v=0$, 
 the Riemann invariants of the observer coincide with the Riemann invariants of the observed system up to the perturbations. 
In contrast, for $\mu^v=1$ 
the observer 
${\bf (R)}$ 
satisfies the same coupling conditions as 
${\bf (S)}$  and no measurement information is inserted.

To investigate the exponential synchronization, we 
consider the difference 
\[\delta := R - S\]
between the state $R$ that is
generated by the observer ${\bf (R)} $ and the original state $S$. 
For the difference  $\delta$ we  have the system
\begin{equation}\label{Diff}\tag{\textbf{Diff}}
\hspace{-.2cm}\left\{
\begin{array}{l}
\delta_\iota^e(0, x)     =  z_\iota^e (x) -   y_\iota^e (x),
\iota \in \{+, \, -,  \, 0 \}, \,  x\in (0,L^e),  \, e\in E,
\\
\\
\delta_{out}^e(t, x^e(v))   =   (1-\mu^v)  Z_{out}^e(t, x^e(v))  + \mu^v \, \delta_{in}^e(t, x^e(v)),
\\
 t\in (0,T),\;   \mbox{\rm if } \;e\in E_0(v)   \mbox{  and  }   |E_0(v)|=1,
  \\
  \delta_0^e(t, x^e(v))   =   Z_{0}^e(t, x^e(v)), 
 \\
 t\in (0,T),\,
 \mbox{\rm if } 
 e\in E_0(v), \; |E_0(v)|=1
 \mbox{  and  }
e\not \in  E_{in}(v,t);
\\ 
\\
\delta_{out}^e(t, x^e(v))     =
  \mu^v\, K^{v,e}_\sigma(\delta_+, \delta_-,t) 
  +
  (1 - \mu^v)  Z_{out}^e(t, x^e(v)) ,
\\
t\in (0,T), 
 \,  \mbox{\rm if } \;
 e\in E_0(v) \mbox{  and  } 
  |E_0(v)|\geq 2,
\\
\delta_0^e(t,\, x^e(v))  = ( 1 - \mu^v)  Z_0^v(t)

\\
+
\mu^v \, [K^v_0(S_0+\delta_0, \, S_+ + \delta_+, \, S_-+\delta_-,\, t)
-
 K^v_0(S_0, \, S_+ , \, S_-,\, t)]
  \\ 
t\in (0,T),  \, \mbox{\rm if }
  e\in E_0(v), \,  |E_0(v)| \geq 2 \,  \mbox{  and  }\, e \not \in E_{in}(v, t) 
\\
\\
 \partial_t
\left(
\begin{array}{r}
\delta_+^e
\\
\delta_-^e
\\
\delta^e_0
\end{array}
\right)
 +
   \Delta^e(S + \delta) \, 
\partial_x
\left(
\begin{array}{r}
\delta_+^e
\\
\delta_-^e
\\
\delta^e_0
\end{array}
\right)
+
{
\color{black}
 \left[  \Delta^e(S + \delta) -   \Delta^e(S) \right]
 }
 \,
 \partial_x 
 \left(
\begin{array}{r}
S_+^e
\\
S_-^e
\\
S^e_0
\end{array}
\right)
\\
 =
\left[
\sigma^e (\delta_+^e+S_+^e,\,\delta_-^e + S_-^e)
-
\sigma^e (S_+^e,\,  S_-^e)
\right]\, 
\left(
\begin{array}{r}
-1
\\
1
\\
0
\end{array}
\right)
\\
\mbox{on }\; [0,T]\times [0,L^e],\, e\in E.
\end{array}
\right.
\end{equation}




\section{A well--posedness result}
\label{wellposedness}

Our system allows to apply a result about the existence of
classical semi-global solutions for quasilinear hyperbolic
systems as stated in
\cite{ci85}, \cite{libook}.
The solutions are constructed 
using integral equations along the characteristic curves
whose slopes are given by the eigenvalues.
A similar  analysis of
Lipschitz solutions for the 
$2\times 2$ system is given in
\cite{gugatulbrich}.
These results 
do not only provide the existence
of a classical solution
on a given finite time interval $[0, \, T]$
for sufficiently small initial
data 
that are $C^1$-compatible
with the boundary conditions 
but also give a  priori bounds
for the $C^1$ norm of the
solutions that are proportional
to the $C^1$-norms of the
initial data.
%
Solutions  of this type have been 
studied by Ta-Tsien Li and his group
in depth, see for example  \cite{lisaintvenant}.

In the quasilinear model that we consider, 
the  eigenvalues in the diagonal system matrix
define three  families of
characteristics.
But it is a particular feature
 of the system that all the eigenvalues only depend on $R_+$ and $R_-$.
This allows to assert  the well-posedness for the
$2\times 2$ system \eqref{pqdiag} 
since the node conditions for $R_+$, $R_-$ are also independent of
$R_0$.
The proof for  the
$2\times 2$ systems with the node conditions from
 ${\bf (S)} $,
  ${\bf (R)} $ and
   ${\bf (Diff)} $
   is very similar to the analysis presented in
   \cite{gugatulbrich}. In
   fact, the only difference  
   are the particular node conditions 
   stated in  ${\bf (R)} $ and
   ${\bf (Diff)} $.
   Note however that the construction in 
 \cite{gugatulbrich} 
 does not yield classical  solutions
but less regular solutions that are Lipschitz continuous with respect to space.
Due to the results from \cite{libook},
we also obtain the existence of classical solutions for the
$2\times 2$ system.
This part of the system also
completely determines
the characteristic curves
with slope $\lambda_0$.
The Riemann invariant $R_0$ is
constant along  these curves.

Under suitable assumptions  our initial data
generate $R_+, R_-$ such that
$|\lambda_0| >0$
everywhere. Then  the sets $E_{in}(v,\, t)$ do not
depend on $t$ 
and the node conditions for $R_0$
uniquely determine the values along the characteristic $\lambda_0$-curves.
The following result is stated in
a $C^1$-neighbourhood of a steady reference state
that is a classical solution of 
 ${\bf (S) }$.
 Steady states in gas networks are discussed in detail in
 \cite{zbMATH06912453}, \cite{nhmneu}.

\begin{theorem}
\label{existence}
Let $T>0$,
a steady state  $(J^e(x))_{e\in E}$
that is a classical solution of 
 ${\bf (S) }$ 
and a 
number $M>0$ be given.
Assume that for all $e\in E$ we have
\begin{equation}
    \label{subsonic}
\lambda_+(J^e_+,\, J^e_-)>0,
\,
\lambda_-(J^e_+,\, J^e_-)<0
\end{equation}
and 
\begin{equation}
    \label{velocityboundinitial}
|J^e_{+}-J^e_-| >0,
\end{equation}
i.e. there is no point on the network with velocity zero. 
Then there exists a number $\varepsilon(T, \, M)>0$ such that
for initial data
$y^e_{\iota}\in C^{1}(0,\, L^e)$
($\iota \in \{+, \,   -, \, 0 \}$,  $e\in E$)
such that
\[
\|y^e_{\pm}-J^e_\pm\|_{ C^1(0,\, L^e) } \leq \varepsilon(T,\, M),\,
\|y^e_{0}- J^e_0\|_{ C^1(0,\, L^e) } \leq \varepsilon(T,\, M),
\]
and control functions
$u^v_\sigma$, $u^e_0 \in C^1(0,\, T)$  ($e\in E$)
such that

\begin{equation}
\label{usigmaassumption}
\left\|
(1- \mu^v)u_{\sigma}^v - {  J_{out}^e(x^e(v))  + \mu^v J_{in}^e(x^e(v) ) }
\right\|_{C^1(0, T)} 
\leq \varepsilon(T,\, M), 
\end{equation}
\begin{equation}
\label{uvassumption}
\|u^e_0
 -    J^e_0(x^e(v)) 
 \|_{C^1(0, T)} \leq \varepsilon(T,\, M)
\end{equation}
that satisfy the  $C^1$-compatibility conditions for 
${\bf (S) }$ 
at the vertices $v\in V$ 
for $t=0$ 
there exists a unique classical solution of
 ${\bf (S) }$
 such that for all $e\in E$ we have
\begin{equation}
\label{lunendlichboundpm}
\|S^e_{\pm} - J^e_\pm\|_{C^1( (0,\, T) \times (0,\, L^e) ) } \leq M,\,
\end{equation}
\begin{equation}
\label{lunendlichbound0}
\|S^e_0- J^e_0\|_{C^1( (0,\, T) \times (0,\, L^e) ) } \leq M
\end{equation}
\begin{equation}
    \label{velocityboundstate}
|S^e_{+}-S^e_-| >0.
\end{equation}

%


\end{theorem}
\textbf{Proof.} 
In order to prove the result,
we consider new variables
$\alpha^e_\pm  =  y^e_{\pm}-J$.
Then  $\alpha^e_\pm  = 0$, $\alpha^e_0  = 0$ 
yields a constant stationary solution on
the network that corresponds to a physical state
with zero flow, constant pressure and no hydrogen.
Moreover,
for $\alpha=0$ the source term in the PDE in
${\bf (S) }$ 
vanishes.

 We normalize the lengths $L^e$ of all the pipes to $1$
without changing the signs of the eigenvalues.
Then the state on each pipe corresponds 
to three components on the interval 
$[0, 1]$.
Thus the state is given by
a function 
with 
${3|E|}$ components  defined on
$[0, 1]$.
After this transformation
all the node conditions on the graph $G$
reduce to boundary conditions
on the transformed
space interval $[0,\, 1]$.
If  $\varepsilon(T, \, M)>0$
is sufficiently small,
we have 
$\lambda_+(S_+^e,S_-^e)>0$
and
$\lambda_+(S^e_+,S^e_-)<0$.
Hence the boundary conditions
for $S_+$ and $S_{-}$ give
the values for the incoming
Riemann invariants as functions
of the outgoing Riemann invariants
with an additive $C^1$-term. 
Hence we can apply Theorem 2.5 in \cite{libook}
(see also Lemma 1.1. in \cite{lirao})
which yield the existence of a unique semi-global
$C^1$-solution of the $2\times 2$ system
\eqref{pqdiag}
that satisfies
the a-priori bound
\eqref{lunendlichboundpm}
for $S_+$ and $S_-$.

By choosing 
$  \varepsilon(T,\, M)>0$ sufficiently small
the a priori bound for $S_+$ and $S_-$
implies that 
we can even obtain a solution such that
\[
|S^e_+ - J^e_+| + |S^e_- - J^e_-| < |J^e_+ - J^e_-|.
\]
Hence due to \eqref{velocityboundinitial}
we can obtain a solution of
the $2\times 2$ system
\eqref{pqdiag}
that satisfies
\[  |S^e_+ - S^e_-| = |S^e_+ - J^e_+ + J^e_+ - J^e_- + J^e_- - S^e_-|
     > 0
\]
and thus \eqref{velocityboundstate}.
Then 
we have
$|\lambda_0(S_+,\, S_-)|>0$.
This implies that for all $t\in [0, T]$ we have
$E_{in}(v,t) = E_{in}(v, \, 0)$.
Hence due to our construction the boundary conditions
for $S_0$  give the necessary values
for the Riemann invariant $R_0$.
%
%
Hence we can apply Theorem 2.5 in \cite{libook}
(see also Lemma 1.1. in \cite{lirao})
to the $3\times 3$ system 
which yield the existence of a unique semi-global
$C^1$-solution that satisfies
\eqref{lunendlichbound0}.



\begin{remark}
\label{existenceremark}
We have discussed 
the well--posedness 
of ${\bf (S) }$.
A similar result can be stated for the
observer system ${\bf (R) }$
with the additional assumption
on the  functions
$Z^e_\iota(t, x_e(v))\in C^1(0, T)$
($v\in V$, $e\in E_0(v)$,
$\iota \in \{+, \,   -, \, 0 \}$)
that should be added to
$u^v_\sigma$, $u^v_0$
in
\eqref{usigmaassumption} and
\eqref{uvassumption}
and should also satisfy the $C^1$-compatibility conditions at
the nodes, i.e. the $Z_\iota^e$ and their derivatives need to vanish for $t \searrow 0$.
These two results also imply the existence of a 
unique classical solution of the error system
${\bf (Diff) }$.

In \cite{gugatulbrich}
a result about the existence of solutions with
the slightly weaker $ W^{1,\infty}$ regularity is shown.
The existence  result for solutions with 
higher $H^2$-regularity
(to be  precise in the space 
$\times_{e\in E} C([0,\, T], \, H^{2}(0, L^e)) $)
for quasilinear systems
given in \cite{bastincoron}
requires $C^2$-regularity of the source term.
In our case, the source term only has $C^1$-regularity.
Therefore, the result from \cite{bastincoron} cannot be applied.

\end{remark}





\section{An $L^2$--Lyapunov function for 
$(\delta_+, \delta_-)$ on the  network}
\label{l2lyapunov}
Squared $L^2$-norms  with exponential weights 
have been used successfully as Lyapunov functions
to show  exponential  boundary feedback stabilization, see
for example \cite{bastincoron}.
We will use a  Lyapunov function of this type 
to show the exponential  synchronization
of the observer system.  
As a first step, we will show the exponential synchronization
of $(R_+, R_-)$ and $(S_+, S_-)$
that are governed by the corresponding 
decoupled $2\times 2$ system.
Later, in a second step we will also show the
exponential synchronization of $R_0$ and $S_0$
with an appropriate Lyapunov function.

For a given real parameter $\psi \geq 0$
and $e\in E$ 
we introduce the weighting functions
\begin{equation}
    \label{hpmdefinition}
h^e_+( x) = B^e_+ \, \exp( -\psi \,x  ), \;\; 
%
h^e_-( x) = B^e_- \,  \exp( \; \; \psi \, x )
\end{equation}
where $ B^e_+$,  $ B^e_-$ are real numbers in $(0,\, \infty)$. 
For $e\in E$,  define the Lyapunov function
\[
{\cal E}^e_\sigma(t) :=  \int_0^{L^e}
h^e_+(x) \, |\delta^e_+(t,\, x)|^2 + 
h^e_-(x) \, |\delta^e_-(t,\, x)|^2 \, dx.
\]
We have
\[
{\cal E}^e_\sigma(t) 
=
\int_0^{L^e}
\begin{pmatrix}
 h^e_+(x) 
\\
 h^e_-(x) 
\end{pmatrix}^\top
\begin{pmatrix}
|\delta^e_+(t,\, x)|^2
\\
 |\delta^e_-(t,\, x)|^2 
\end{pmatrix}   \, dx. 
\]

Note that
${\cal E}^e_\sigma(t)$ is equivalent to the squared $L^2$-norm of $\delta^e$.
Assume that on the time-interval $[0, \, T]$,
a semi-global classical solution of
${\bf (Diff) }$ is given.
Then for the time-derivative of  ${\cal E}^e_\sigma(t)$  we have
\[
\frac{d}{dt} {\cal E}^e_\sigma(t) 
=
2
\int_0^{L^e}
\begin{pmatrix}
 h^e_+(x) 
\\
 h^e_-(x) 
\end{pmatrix}^\top
\begin{pmatrix}
  \delta^e_+(t,\, x) \, \partial_t  \delta^e_+(t,\, x) 
\\
  \delta^e_-(t,\, x) \, \partial_t  \delta^e_-(t,\, x) 
\end{pmatrix}   \, dx. 
\]

Using the partial  differential  equation from ${\bf (Diff) }$,
we will show  in the proof of Theorem \ref{l2lyapunovthm} 
 that the network Lyapunov function
\[
{\cal E}_\sigma(t) := 
\sum_{e\in E} {\cal E}^e_\sigma(t)
\]
satisfies a differential inequality of the form
\[
\frac{d}{dt}
{\cal E}_\sigma(t)
\leq
-  \chi \,{\cal E}_\sigma(t)
- Q(t)
\]
with a positive real number $\chi >0$
and 
\begin{equation}
\label{qdefinition}
Q(t)=
\sum_{e\in E}
\left[
  h^e_+(x) \, \lambda_+^e( S +  \delta ) \,  |\delta^e_+|^2(t,\, x) 
+ h^e_-(x) \, \lambda_-^e( S +  \delta ) \,  |\delta^e_-|^2(t,\, x)
\right]|_{x=0}^{L^e}.
\end{equation}
If we have an inequality of the form
 $-Q(t)\leq \eta$
 with a real number $\eta$,
 this yields the exponential decay of
${\cal E}_\sigma(t)$ 
up to the level of perturbations with the following version of Gronwall's Lemma:

\begin{lemma}
\label{gw1}
Let $\chi>0$, $\eta\geq 0$
and a function ${\cal E}(t) $ that is absolutely contiuous on
$[0, T]$ be given.
Assume that for $t\in [0, T]$ 
almost everywhere we have 
\[
{\cal E}'(t)
\leq
- \chi \, {\cal E}(t) + \eta.
\]
Then for all $t\in [0, T]$ we have 
\begin{equation}
    \label{lpfctnungl}
{\cal E}(t) 
\leq
{\cal E}(0) \, \exp(  - \chi  \, t)
+
 \frac{\eta}{ \chi  } 
\left[ 1 -\exp(  - \chi  \, t)
\right]
.
\end{equation}
\end{lemma}
\begin{proof}
Define
${\cal H}(t) =  \exp( \chi \, t)\, \left[  {\cal E}(t) - \frac{\eta}{ \chi  }  \right] $. The product rule yields
\begin{eqnarray*}
{\cal H}'(t) 
& =   &
\chi  \, \exp(  \chi \, t)\,  \left[  {\cal E}(t) - \frac{\eta}{ \chi  }   \right] 
+ \exp(  \chi\, t) \,   {\cal E}'(t)  
\\
& \leq &
\exp(  \chi  \, t)
\left[
  \chi  \,  \left[  {\cal E}(t) - \frac{\eta}{  \chi   }   \right] 
- \chi \, {\cal E}(t) + \eta
\right]
\\
& = &
0
\end{eqnarray*}
for $t$ almost everywhere in $[0,\, T]$. 
We get 
$
{\cal H}(t) \leq 
{\cal H}(0)
.
$
Hence due to the definition of $ {\cal H} $ we 
obtain
$
{\cal E}(t)  -  \tfrac{\eta}{ \mu    } 
\leq
\left[ {\cal E}(0) -  \tfrac{\eta}{ \mu    } \right]  \, \exp(  - \chi  \, t)
.
$
Thus we have shown Lemma \ref{gw1}. 
%
\end{proof}
%
%
In order to study the properties of $Q(t)$ we
reorder the sum in terms of the vertices of the graph and obtain 
\begin{multline}
Q(t)=
\sum_{v\in V}
\sum_{e\in E_0(v): x^e(v) = L^e}
\begin{pmatrix}
 h^e_+(L^e) 
\\
 h^e_-(L^e) 
\end{pmatrix}^\top
\begin{pmatrix}
 \lambda_+^e( S +  \delta )  |\delta^e_+|^2(t,\, L^e)|
\\
 \lambda_-^e( S +  \delta )  |\delta^e_-|^2(t,\, L^e)|
\end{pmatrix}  
\\
-
\sum_{e\in E_0(v): x^e(v) = 0}
\begin{pmatrix}
 h^e_+(0) 
\\
 h^e_-(0) 
\end{pmatrix}^\top
\begin{pmatrix}
 \lambda_+^e( S +  \delta )  |\delta^e_+|^2(t,\, 0)|
\\
 \lambda_-^e( S +  \delta )  |\delta^e_-|^2(t,\, 0)|
\end{pmatrix}.
\end{multline}
Define
$I_{L}(v) = \{e\in E_0(v): x^e(v)=L^e\}$ and
$I_{0}(v) = \{e\in E_0(v): x^e(v)=0\}$.
Using the term for the nodal coupling in  the definition of  ${\bf (D{iff})}$
with the parameter $\mu^v$ this yields 

\begin{multline}
Q(t) =
\sum_{v\in V}\Bigg[
\sum_{e\in I_{L}(v)}
\Big[
  h^e_+( L^e) \, \lambda_+^e( S +  \delta )(t,\, L^e) \,  |\delta^e_{in}(t,\, L^e)|^2
  \\
  +h^e_-(L^e) \, \lambda_-^e( S +  \delta )(t,\, L^e)   \,  
  \left| \mu^v 
 K_\sigma^{v,e}(\delta_+,\delta_-,t) 
 + (1 - \mu^v) Z^e_{out}(t, \, L^e) \right|^2
 \Big]
\\
-
\sum_{e\in I_{0}(v)}
\Big[
  h^e_+( 0 ) \, \lambda_+^e( S +  \delta )(t,\, 0)   \, 
 \left| \mu^v 
 K_\sigma^{v,e}(\delta_+,\delta_-,t) 
 + (1 - \mu^v) Z^e_{out}(t, \, 0) \right|^2
\\
+
h^e_-( 0 ) \, \lambda_-^e( S +  \delta )(t,\, 0) \,  |\delta^e_{in}(t,\, 0)|^2\Big].\Bigg]
\end{multline}

We define 
\begin{eqnarray*}
\underline C^v(t)&:=
\min\{ 
\min_{e\in I_L(v)} h^e_+( L^e) \, \lambda_+^e( S +  \delta )(t,\, L^e), \,
\min_{e\in I_0(v)} h^e_-( 0 ) \, |\lambda_-^e( S +  \delta )(t,\, 0)|
\}\\
\overline C^v(t)&:=
\max\{ 
\max_{e\in I_L(v)} h^e_-( L^e) \, |\lambda_-^e( S +  \delta )(t,\, L^e)|, \,
\max_{e\in I_0(v)} h^e_+( 0 ) \, \lambda_+^e( S +  \delta )(t,\, 0)
\}.
\end{eqnarray*}
and observe that $\underline  C^v(t) >0 $. 
Note that
$K_\sigma^{v,e}(\cdot,\,\cdot,t)$ is a linear map and 
there exist a constant $\hat C^v$ such that
\[ 
\|K_\sigma^{v,e}(\delta_+,\delta_-,t) \| \leq 
\hat C^v \, \|\delta_{in}^v(t)\|.\]
Then we have
\[Q(t)
\geq
\sum_{v\in V}
\left[\underline C^v(t)
- 2 \,\overline C^v(t) \, \hat C^v |\mu^v|^2
\right]\|\delta_{in}^v(t)\|^2
- 2 | 1 - \mu^v|^2 \,\|Z_{out}^v(t)\|^2.
\]

If for a given value of $\psi$
(which  influences the values of
$ h^e_\pm( 0 ) $ and
$ h^e_\pm( L^e ) $) 
the nodal parameter $\mu^v$ is chosen in 
such a way that $ |\mu^v |$ is sufficiently small, we have 
\begin{equation}
    \label{inequality2023}
\underline C^v(t)
- 2 \,\overline C^v(t) \, \hat C^v |\mu^v|^2
\geq 0
 \end{equation}
and thus
\begin{equation}
\label{qbound}
Q(t)
\geq 
\sum_{v\in V}
- 2 | 1 - \mu^v|^2 \,\|Z_{out}^v(t)\|^2.
\end{equation}
For all $e\in E$, define the number
\begin{equation}
\label{kappadefinition}
\kappa^e ;=  \max_{z \in [0, \, L^e]} 
\left\{\frac{h^e_+(z)}{h^e_-(z)},
\,
\frac{h^e_-(z)}{h^e_+(z)}
\right\} \geq 1.
\end{equation} 
The following Theorem contains
our result about the exponential synchronization
of $(R_+, R_-)$ and $(S_+, S_-)$.
Similar to Theorem \ref{existence} it
is stated in some
neighbourhood of a steady reference state
that is a classical solution of the first two components of
 ${\bf (S) }$.
 The results in
 \cite{zbMATH06912453}, \cite{nhmneu}
 imply that for sufficiently small gas velocities
 the norms of the derivative
 in the steady state become arbitrarily small.
 To be precise, we can choose the  steady state in 
such a way that 
$(J^e_+, J^e_-)$
is arbitrarily close to a constant function
of the form 
$(C, C)$
for all $e\in E$.


%
\begin{theorem}
\label{l2lyapunovthm}
Let 
$M>0$ and
a classical steady state $(J^e_+,\, J^e_-)_{e\in E}$ 
 of ${\bf (S)}$ 
that satisfies the conditions in Theorem \ref{existence} 
be given.
Assume that for all $e\in E$ we have 
\begin{equation}
    \label{jassumption}
    |J^e_+ - J^e_-| \leq M.
\end{equation}
Let $T > 0$  be given.
For all $e\in E$, let initial states 
$y^e_+$,  $y^e_-$, $z^e_+$, $z^e_-
\in  C^1(0, \, L^e)$ be given
such that at all nodes $ v\in V$  the $C^1$-compatibility conditions
are satisfied.
Assume that the first two components  of  ${\bf (S) }$ and ${\bf (R)}$
have classical solutions on $[0, T]$
such that 
(\ref{lunendlichboundpm})
and the following
a-priori bounds
\eqref{aprioribound2a}, 
\eqref{aprioribound1} 
and 
\eqref{aprioribound2}
hold:

%
%




There exists  numbers 
$\tMz>0$, 
$c>0$ and
$\varepsilon_0>0$
such that
%
\begin{equation}
\label{aprioribound2a}
|\partial_x  S^e_\pm| \leq \tMz,
\end{equation}
%
\begin{equation}
\label{aprioribound1}
  \frac{3}{4}\,c \leq  \lambda_+^e( R) 
  \leq \frac{5}{4} \, c
  ,
\;\;
   -\frac{5}{4}\,c \leq  
 \lambda_-^e( R) \leq -\frac{3}{4} \, c,
\end{equation}


\begin{equation}
\label{aprioribound2}
 |\partial_x \left( \lambda_\pm^e( R ) \right)| \leq \varepsilon_0
. 
\end{equation}

Assume  that we have a nonempty compact  convex set $ {\cal U} \subset {\mathbb R^2}$
that contains the function values
of the solutions 
 $S^e$ and $R^e$ and a number $\beta>0$ such that
\begin{equation}
\label{aprioribound4}
\max_{ \tilde R \in {\cal U}  } \left\{ \left|\partial_{R_+} \lambda_\pm^e( \tilde R )\right|,\,
\left|\partial_{R_-} \lambda_\pm^e( \tilde R )\right| \right\}
\leq \beta.
\end{equation}

Then on $[0, \, T]$ there exists a classical
solution $(\delta_+, \delta_-)$
of the first two components 
of system ${\bf (Diff)}$.
%
\\
Let $\psi \in (0, \, \infty)$ 
be given such that  for all $e\in E$ we have 
(with $\kappa^e$ as defined in \eqref{kappadefinition}) 
\begin{equation}
\label{10012023}
\varepsilon_0 +
( 12 \, \nu^e \, 
M
+ 4 \,  \beta \, \tMz )
\, 
 \kappa^e
<    \frac{3}{4}\, \psi \,c.
\end{equation}
Assume that for each  node
$v\in V$  the  number 
$\mu^v\in [0, 1]$ is 
sufficiently small in the sense that
for all $t\in [0,\, T]$
inequality \eqref{inequality2023} holds.
%
%
Assume that the perturbations are bounded from above by 
$\eta_\sigma$, in the sense that for all $t\in [0,\, T]$ we have 
\begin{equation}
    \label{pertbound}
\sum_{v\in V}
 2 \,\|Z_{out}^v(t)\|^2
\leq \eta_\sigma.
\end{equation}
Then  the  solution of  ${\bf (Diff)}$
%
decays   exponentially fast 
up to the perturbation level 
in the sense that
there exist  constants
$C_1 >0$, 
$\chi>0$
such that for all $t\in [0,\, T]$ 
we have
\begin{multline}
\label{shirilyapunov}
\sum_{e\in E}
 \int_0^{L^e} \left|\delta^e_+(t,\, x)  \right|^2 
 + \left| \delta^e_-(t,\, x) \right|^2
\, dx\\
\leq C_1
\exp(-\chi \,  t) \,
\sum_{e\in E}
 \int_0^{L^e} \left|\delta^e_+(0,\, x)  \right|^2
 + \left| \delta^e_-(0,\, x) \right|^2 \, dx
 + C_1\, \frac{\eta_\sigma}{\chi}.
\end{multline}
Hence the $L^2$-norm of the error $\delta$ between the  state $R$
of the observer and the  state  $S$ of the original system  decays exponentially fast
up to the perturbation level.
\end{theorem}

\begin{remark}
We want to emphasize that the inequalities
\eqref{jassumption}
(\ref{aprioribound2a}), 
(\ref{aprioribound1}),
(\ref{aprioribound2}) 
hold 
on the time interval $[0, T]$ 
for all $e\in E$ on the intervals $[0,\, L^e]$.

Inequality 
(\ref{aprioribound2a})
is satisfied if
the steady reference state satisfies 
$|\partial_x  J^e_\pm| \leq \tMz/2$ and
$M\leq \tMz/2$ which yields 
$|\partial_x  (S^e_\pm- J^e_\pm)| \leq \tMz/2$.
Analogously, we can find
$\varepsilon_0>0$ such that
(\ref{aprioribound2}) holds.
Due to \eqref{subsonic}, we can find $c>0$ such that
\eqref{aprioribound1} is satisfied.

Inequality (\ref{10012023}) is satisfied if
$\varepsilon_0$, 
$M$
and $\tMz $ are sufficiently small.
This can be achieved  in a neighbourhood of the steady state
 $(J^e_+, J^e_-)$ if it is sufficiently close to a constant function.
This is the case if the values of $\nu^e>0$ are sufficiently small
or if the flow rates are sufficiently small.
%

For the synchronization result for
$(\delta_+, \, \delta_-)$
it is not necessary that
the steady reference state satisfies
\eqref{velocityboundinitial}.
The assertion of Theorem \ref{l2lyapunovthm}
also holds for steady states of the form
$(J^e_+, J^e_-) = (C, C)$
for all $e\in E$.
This means that the
$2\times 2$ observer system
for $(\delta_+, \, \delta_-)$ can also
determine the direction of the flow.
Equation \eqref{velocityboundinitial} is only necessary  for the 
analysis of the
synchronization of $\delta_0$ that is presented
below.
\end{remark}

\begin{proof}
Since $R:= S + \delta$,
(\ref{aprioribound1}) and
(\ref{aprioribound2})
can also be seen as conditions on 
$\delta$.
Define 
\begin{multline}
F^e(S, \, \delta)
:=
\left(
\sigma^e (\delta_+^e+S_+^e,\,\delta_-^e + S_-^e)
-
\sigma^e (S_+^e,\,  S_-^e)
\right)
\left(
\begin{array}{r}
-1
\\
1 \end{array}
\right)
\\
-
 \left[ \hat \Delta^e(S + \delta) -  \hat \Delta^e(S) \right]
 \left(
\begin{array}{r} 
\partial_x \, S^e_+
\\
\partial_x \, S^e_-
\end{array}
\right)
.
\end{multline}
The first two components of the 
partial  differential  equation in ${\bf (Diff) }$ can be 
stated as 
{\small
\[\partial_t
\begin{pmatrix}
\delta^e_+
\\
\delta^e_-
\end{pmatrix}
+
\begin{pmatrix}
\lambda_+( S^e_+ +  \delta^e_+,  S^e_- +  \delta^e_-) &  0
\\
0 & \lambda_-(  S^e_+ +  \delta^e_+,  S^e_- +  \delta^e_-)  
\end{pmatrix}
\partial_x
\begin{pmatrix}
\delta^e_+
\\
\delta^e_-
\end{pmatrix}
=
\begin{pmatrix}
F^e_+( S, \,  \delta)
\\
F^e_-( S, \,  \delta)
\end{pmatrix}.
\]
} 
Since the PDE is in diagonal form this yields
\begin{multline}
\frac{d}{dt}{\cal E}^e_\sigma(t)  
=
2 \int_0^{L^e}
\begin{pmatrix}
 h^e_+(x)  
\\
 h^e_-(x) 
\end{pmatrix}^\top
\left[
-
\begin{pmatrix}
  \lambda_+^e( S +  \delta )\, \delta^e_+ \, \partial_x \delta^e_+
\\
  \lambda_-^e( S +  \delta ) \, \delta^e_- \, \partial_x \delta^e_-
\end{pmatrix} 
+
\begin{pmatrix}
\delta^e_+ \, F^e_+( S, \,  \delta)
\\
 \delta^e_- \, F^e_-( S, \,  \delta)
\end{pmatrix}
\right]\, dx\\
=
 \int_0^{L^e}
\begin{pmatrix}
 h^e_+(x) 
\\
 h^e_-(x) 
\end{pmatrix}^\top
\left[
-
\begin{pmatrix}
  \lambda_+^e( S +  \delta ) \, \partial_x (\delta^e_+)^2
\\
  \lambda_-^e( S +  \delta ) \, \partial_x (\delta^e_-)^2
\end{pmatrix} 
+
2
\begin{pmatrix}
\delta^e_+ \,F^e_+( S, \,  \delta)
\\
\delta^e_- \,F^e_-( S, \,  \delta)
\end{pmatrix}
\right]\, dx.
\end{multline}


Using  integration by parts we obtain 
\begin{multline}
\frac{d}{dt}{\cal E}^e_\sigma(t)  
=
 \int_0^{L^e}
\begin{pmatrix}
\partial_x \left[ h^e_+(x) \, \lambda_+^e( S +  \delta ) \right]
\\
\partial_x \left[ h^e_-(x) \, \lambda_-^e( S +  \delta ) \right]
\end{pmatrix}^\top
\begin{pmatrix}
  (\delta^e_+)^2
\\
  (\delta^e_-)^2
\end{pmatrix} 
\,dx
\\
-
\left.
\begin{pmatrix}
 h^e_+(x) \, \lambda_+^e( S +  \delta )  
\\
h^e_-(x) \, \lambda_-^e( S +  \delta )  
\end{pmatrix}^\top 
\begin{pmatrix}
  (\delta^e_+)^2
\\
  (\delta^e_-)^2
\end{pmatrix} 
\right|_{x=0}^{L^e}
+
2
 \int_0^{L^e}
\begin{pmatrix}
 h^e_+(x) 
\\
 h^e_-(x) 
\end{pmatrix}^\top
\begin{pmatrix}
\delta^e_+ \, F^e_+( S, \,  \delta)
\\
\delta^e_- \,F^e_-( S, \,  \delta)
\end{pmatrix}
\, dx.
\end{multline}

We have
\[\partial_x \left( h^e_\pm(x) \, \lambda_\pm^e( S +  \delta ) \right)=
h^e_\pm(x) \left( \mp \psi \, \lambda_\pm^e( S +  \delta )  +
\partial_x (\lambda_\pm^e( S +  \delta ))
 \right).
\]

Due to  (\ref{aprioribound1})  and  (\ref{aprioribound2}) we have
\[
\partial_x \left( h^e_\pm(x) \, \lambda_\pm^e( S +  \delta ) \right)
\leq 
- h^e_\pm(x) \left( \psi \,\, \frac{3}{4} \, c  -
 \varepsilon_0
 \right).
\]
This yields
the inequality
\begin{equation}
\label{01122022}
\int_0^{L^e}
\partial_x 
\begin{pmatrix}
 h^e_+(x) \, \lambda_+^e( S +  \delta ) 
\\
 h^e_-(x) \, \lambda_-^e( S +  \delta ) 
\end{pmatrix}^\top
\,
\begin{pmatrix}
|\delta^e_+|^2
\\
 |\delta^e_-|^2 
\end{pmatrix}   \, dx
\leq
-  \left( \frac{3}{4}  \psi \, c  - \varepsilon_0   \right) \, {\cal E}^e_\sigma(t).
\end{equation}

Since the right-hand side of the inequality is negative, it is
a suitable start to derive a differential inequality
that can be used to apply Gronwall's lemma.
It remains to derive bounds for the other terms that appear in
$\frac{d}{dt}{\cal E}^e_\sigma(t)$.
We have
\begin{multline}
  2  \int_0^{L^e}
  \begin{pmatrix}
h^e_+(x) 
\\
h^e_-(x)
\end{pmatrix}^\top 
\begin{pmatrix}
\delta^e_+ \, F^e_+( S, \,  \delta)
\\
\delta^e_- \, F^e_-( S, \,  \delta)
\end{pmatrix}
\, dx
\\
 = 2  \int_0^{L^e}
\left[\sigma^e (\delta_+^e+S_+^e,\,\delta_-^e + S_-^e)
-
\sigma^e (S_+^e,\,  S_-^e)\right] 
   \begin{pmatrix}
h^e_+(x) 
\\
h^e_-(x)
\end{pmatrix}^\top 
\begin{pmatrix}
-\delta^e_+ 
\\
\delta^e_- 
\end{pmatrix}
\, dx
\\
- 
  2  \int_0^{L^e}
  \begin{pmatrix}
h^e_+(x) \, \delta^e_+
\\
h^e_-(x) \, \delta^e_-
\end{pmatrix}^\top 
\left[
  \hat\Delta^e(S + \delta) -   \hat\Delta^e(S) \right]
  \,
\left(
\begin{array}{r} 
\partial_x \, S^e_+
\\
\partial_x \, S^e_-
\end{array}
\right)
\, dx.
\end{multline}
Let us define the term that is  connected with  the 
source term that models the friction at the pipe walls
\[
I_1^e :=
 2  \int_0^{L^e}
\left[\sigma^e (\delta_+^e+S_+^e,\,\delta_-^e + S_-^e)
-
\sigma^e (S_+^e,\,  S_-^e)\right] 
   \begin{pmatrix}
h^e_+(x) 
\\
h^e_-(x)
\end{pmatrix}^\top 
\begin{pmatrix}
-\delta^e_+ 
\\
\delta^e_- 
\end{pmatrix}
\, dx.
\]

Moreover, we define
\[I_2^e
:=
- 
  2  \int_0^{L^e}
  \begin{pmatrix}
h^e_+(x) \, \delta^e_+
\\
h^e_-(x) \, \delta^e_-
\end{pmatrix}^\top 
\left[
  \hat\Delta^e(S + \delta) -   \hat\Delta^e(S) \right]
\left(
\begin{array}{r} 
\partial_x \, S^e_+
\\
\partial_x \, S^e_-
\end{array}
\right)
\, dx.
\]

Now we derive an upper bound for
$|I_1^e|$. 
Since the function $z\mapsto \nu^e \, z\, |z|$ is
differentiable with the derivative
$z \mapsto 2\, \nu^e \, |z|$, 
definition 
(\ref{sigmadefinition})
of the function $\sigma^e$,
the mean value theorem,
\eqref{jassumption}
and
\eqref{lunendlichboundpm}
yield
\[
\left|\sigma^e (\delta_+^e+S_+^e,\,\delta_-^e + S_-^e)
-
\sigma^e (S_+^e,\,  S_-^e)\right|
\leq
2\, \nu^e \,
(3 M) 
 |\delta_+^e - \delta_-^e|.
\]

Hence we have 
\begin{multline}
|I_1^e | \leq 6 \, \nu^e \,M
\int_0^{L^e}
\left[
h^e_+(x) \, \left| \delta^e_+\right|
+
 h^e_-(x) \, \left| \delta^e_-\right|
\right]
\left| \delta^e_+ - \delta^e_-\right|
\, dx
\\
\leq 6 \, \nu^e \,M
\,
\int_0^{L^e}
\left[
h^e_+(x) \,  \left| \delta^e_+\right|
+
h^e_-(x) \, \left| \delta^e_-\right|
\right]
\left[
\left| \delta^e_+\right|  + \left|\delta^e_-\right|
\right]
\,dx
\\
= 6 \, \nu^e \,
M 
\int_0^{L^e}
\left[
h^e_+(x) \,  \left| \delta^e_+\right|^2
+
h^e_-(x) \,  \left| \delta^e_-\right|^2
\right]
+ [h^e_+(x) + h^e_-(x) ]\left| \delta^e_+ \delta^e_-\right|
\,dx
\\
\leq 6 \, \nu^e \,
M
\int_0^{L^e}
\left[\frac{3}{2} h^e_+(x) + \frac{1}{2} h^e_-(x)\right]  \left| \delta^e_+\right|^2
+
\left[\frac{1}{2} h^e_+(x) + \frac{3}{2} h^e_-(x)\right] 
 \left| \delta^e_-\right|^2
\,dx.
\end{multline}
Since 
$h^e_\pm(x) =  \frac{h^e_\pm(x)}{h^e_\mp(x)} \, h^e_\mp(x)$,
we have
$h^e_\pm(x) \leq \kappa^e \,
\, h^e_\mp(x)$. 
Thus we have the bound 
$$
|I_1^e |
\leq 6\, 
\nu^e \, M
2\,  \kappa^e 
\, 
{\cal E}^e_\sigma(t).
$$

It remains to derive a bound for the term $|I_2^e| $ that comes from the variations
of the eigenvalues of the system matrix.
Due to 
(\ref{aprioribound2a})
and
(\ref{aprioribound4}) 
we have
\begin{multline}
|I_2^e| \leq 
 2  
 \int_0^{L^e}
[ h^e_+(x) \, |\delta^e_+| + h^e_-(x) \, |\delta^e_-|]
\, 
\beta \, [|\delta^e_+| + |\delta^e_-| ]\tMz 
\, dx
\\
=
 2 \,  \beta \,  \tMz 
 \int_0^{L^e}
 h^e_+(x) \, |\delta^e_+|^2 + h^e_-(x) \, |\delta^e_-|^2
+ [h^e_+(x) + h^e_-(x) ] \,|\delta^e_+\, \delta^e_-| ]
\, dx.
\end{multline}

Similarly as for $I_1^e$ we obtain the inequality
\begin{multline}
|I_2^e|  
\leq 2 \,  \beta \,  \tMz \, 
\int_0^{L^e}
\left[\frac{3}{2} h^e_+(x) + \frac{1}{2} h^e_-(x)\right]  \left| \delta^e_+\right|^2
\\+
\left[\frac{1}{2} h^e_+(x) + \frac{3}{2} h^e_-(x)\right] 
 \left| \delta^e_-\right|^2
\,dx
\leq
 2 \,  \beta \,  \tMz \, 
 2\,  \kappa^e 
\, 
{\cal E}^e_\sigma(t)
 .
\end{multline}

We have
$
  2  \int_0^{L^e}
  \begin{pmatrix}
h^e_+(x) 
\\
h^e_-(x) 
\end{pmatrix}^\top 
\begin{pmatrix}
\delta^e_+ \, F^e_+( S, \,  \delta)
\\
\delta^e_- \, F^e_-( S, \,  \delta)
\end{pmatrix}
\, dx
= I_1^e + I_2^e.
$

Hence 
\eqref{01122022} and the bounds for $|I_1^e|$ and $|I_2^e|$
yield
the differential inequality
\begin{multline*}
\frac{d}{dt}{\cal E}^e_\sigma(t)  
\leq
\left[
-  \frac{3}{4}\, \psi \,c  +   \varepsilon_0 
+
( 6 \, \nu^e \, 
M
+ 2 \,  \beta \, \tMz )
2\,  \kappa^e 
 \right]
\, {\cal E}^e_\sigma(t) 
\\
-
\left[
  h^e_+(x) \, \lambda_+^e( S +  \delta ) \,  |\delta^e_+(t,\, x)|^2
+ h^e_-(x) \, \lambda_-^e( S +  \delta ) \,  |\delta^e_-(t,\, x)|^2
\right]|_{x=0}^{L^e}.
\end{multline*}
We define the number 
\[
\chi := \min_{e\in E} \frac{3}{4}\, \psi \,c  -   \varepsilon_0 
-
( 6 \, \nu^e \, 
M
+ 2 \,  \beta \, \tMz )
2\,  \kappa^e 
.
\]
Then \eqref{10012023} implies $\chi >0$.
For  the network Lyapunov function we obtain
$$\frac{d}{dt} {\cal E}_\sigma(t) \leq
- \chi 
\,
{\cal E}_\sigma(t) - Q(t)
$$
with $Q(t)$ as defined in \eqref{qdefinition}.
Since  \eqref{inequality2023} holds,
due to \eqref{pertbound} we obtain
$$
\frac{d}{dt} {\cal E}_\sigma(t) \leq
- \chi 
\,
{\cal E}_\sigma(t) + \eta_\sigma.
$$
Now Lemma  \ref{gw1} implies
\eqref{lpfctnungl}.
Since  
${\cal E}_\sigma(t) $
is equivalent to the $L^2$-norm, we obtain
\eqref{shirilyapunov}.

\end{proof}

\section{An $L^2$--Lyapunov function for  $\delta_0$ on the  network}
\label{l2lyapunovh}

In order to complete our analysis for the case of hydrogen blending 
in this section we discuss a Lyapunov function
for $\delta_0$.
Note that in order to show the exponential decay of
$\delta_0$, we 
use the exponential decay of $(\delta_+, \delta_-)$ that
we have shown in the previous section.
For our analysis, we use the following version of Gronwall's Lemma:
\begin{lemma}
\label{gwl2}
Let 
 $\chi_0>0$, $\chi>0$, $\eta_0\geq 0$, $\eta\geq 0$,  $D_0>0$ 
and a function ${\cal E}(t) $ that is absolutely contiuous on
$[0, T]$ be given.
Assume that 
\begin{equation}
    \label{gwl2vor}
{\cal E}'(t)
\leq
- \chi_0 \, {\cal E}(t) 
+ 
D_0 \, 
\left[\sqrt{{\cal E}(t)} \, 
\sqrt{
\exp( - {\chi}\, t)
+
\martin{
 \eta
 }
 }
+
\martin{
\exp( - {\chi}\, t)
+
 \eta
 }
 \right]
+
\eta_0
\end{equation}
for $t\in [0, \, T]$ almost everywhere.
%
%
Then for all 
 $\zeta < \chi_0$,
 $\zeta \not = \chi$ 
and
$t\in [0, T]$ we have 
\begin{equation}
    \label{gw2}
{\cal E}(t) 
\leq
{\cal E}(0) \, 
\textrm{e}^{ - \zeta \, t}
+
 \frac{\eta_0}{ \zeta
 } 
+
\left[
\frac{ D_0^2}{4(\chi_0 - \zeta)  } \, \frac{ 1 }{    (\zeta - \chi)  } 
+
D_0
\right]
\left[\frac{ 
\textrm{e}^{ - \chi \, t} - \textrm{e}^{ -\zeta \, t}
}{    \zeta - \chi  } 
+
\frac{\eta}{\zeta}
\left(
1 -
\textrm{e}^{ -\zeta \, t}
\right)
 \right]
 .
 \end{equation}
\end{lemma}
\begin{proof}
Let $\zeta < \chi_0$ be given. Young's inequality  implies 
$$
\sqrt{{\cal E}(t)} \, 
\frac{  \sqrt{ 2 (\chi_0 - \zeta)}}{\sqrt{2 ({\chi_0 - \zeta})}} \,D_0 \, 
\sqrt{
\exp( - {\chi}\, t)
+
\martin{
 \eta
 }
 }
\leq
{\cal E}(t) (\chi_0 - \zeta) 
+
\frac{ D_0^2 }{ 4 (\chi_0 - \zeta)}  \, 
\left[
\exp( - {\chi}\, t)
+
\martin{
 \eta
 }
 \right]
.
$$
Thus we have
$
{\cal E}'(t)
\leq
 - \zeta\, {\cal E}(t) 
+
\left[ \frac{ D_0^2 }{4(\chi_0 - \zeta)}  
+
D_0\right]\, 
 \left[
\exp( - {\chi}\, t)
+
\martin{
 \eta
 }
 \right]
+
\eta_0
$
.

Define
${\cal H}(t) =  \exp(  \zeta\, t)\, \left[  {\cal E}(t) - \frac{\eta_0}{  \zeta }  \right] $. 
The product rule yields

\begin{eqnarray*}
{\cal H}'(t) 
& =   &
\zeta \, \exp(   \zeta \, t)\,  \left[  {\cal E}(t) - \frac{\eta_0}{  \zeta  }   \right] 
+ \exp(  \zeta \, t) \,   {\cal E}'(t)  
\\
& \leq &
\exp(   \zeta  \, t)
\left[
  \zeta  \,  \left[  {\cal E}(t) - \frac{\eta_0}{  \zeta }   \right] 
- \zeta \, {\cal E}(t) + \eta_0
 +
 \left[\frac{ D_0^2 }{4(\chi_0 - \zeta)}  + D_0 \right] \, 
 \left[
\exp( - {\chi}\, t)
+
\martin{
 \eta
 }
 \right]
\right]
\\
& = &
\left[\frac{ D_0^2 }{{4( \chi_0 - \zeta) }} + D_0 \right]  \, 
\left[ \exp\left( (\zeta - \chi) \, t\right)
+
\eta\,
\exp(\zeta\, t)
\right]
\end{eqnarray*}
for $t\in [0,\, T]$ almost everywhere. 
 This yields 
$$
{\cal H}(t) \leq 
{\cal H}(0)
+
\left[
 \frac{ D_0^2}{4(\chi_0 - \zeta)  } 
 +
 D_0
 \right]\, 
\left[ 
\frac{ \exp(  \left(\zeta - \,\chi \right) t) - 1  }{    \zeta - \chi  } 
+
\frac{\eta}{\zeta}
\left(
\exp(\zeta \, t) - 1
\right)
\right]
.
$$
Hence due to the definition of $ {\cal H} $ we have
$$
{\cal E}(t)  -  \frac{\eta}{ \zeta   } 
\leq
\left[ {\cal E}(0) -  \frac{\eta}{ \zeta   } \right]  \, 
\textrm{e}^{ - \zeta \, t}
+
\left[\frac{ D_0^2}{4(\chi_0 - \zeta)  } + D_0  \right]
\left[\frac{ 
\textrm{e}^{ - \chi \, t} - \textrm{e}^{ -\zeta \, t}
}{    \zeta - \chi  } 
+
\frac{\eta}{\zeta}
\left(
1 -
\textrm{e}^{ -\zeta \, t}
\right)
 \right]
$$
and \eqref{gw2} follows.
Thus we have proved Lemma \ref{gwl2}.

%

\end{proof}

In order to analyze the exponential decay of $\delta_0$,
we assume that
the first two components 
$(\delta_+, \delta_-)$
of the solution of   ${\bf (Diff)} $ are given
and that for all $e\in E$, $x\in [0, \, L^e]$
and $t\in [0,\, T]$ we have
\begin{equation}
    \label{signumvor}
|\lambda_0(S^e_+ + \delta^e_+,\,S^e_- + \delta^e_-)(t, \, x) )|>0.
\end{equation}
Due to 
\eqref{velocityboundstate}
in Theorem \ref{existence} we know that this assumption
holds in a $C^1$-neighborhood of
a steady state that has non-zero velocity everywhere.
We use the standard notation
${sign}
(x) = 1 $ if $x>0$,
${sign}
(x) = - 1$ if $x<0$
and
${sign}
(0) = 0$.
Due to \eqref{signumvor} for $e\in E$
the number 
$$s^e  :=  {sign}(\lambda_0(S^e_+ + \delta^e_+,\,S^e_- + \delta^e_-)(t, \, x) ))$$
is well-defined.

Let  $\psi_0 \geq 0$  denote a real parameter.
For $e\in E$ define
\begin{equation}
\label{e0defi}
{\cal E}^e_0(t)
 =  
 \int_0^{L^e} 
 \textrm{e}^{- s^e \,\psi_0\, x}
 \, |\delta^e_0(t,\, x)|^2\, dx
.
\end{equation}
We will use
${\cal E}^e_0(t)$
as a Lypunov function to show the following theorem
about the synchronization of $\delta_0$:

\begin{theorem}
\label{l2lyapunovthm0}
Let 
$M>0$ and
a classical steady state $(J^e_+,\, J^e_-,J^e_0)_{e\in E}$ 
 of ${\bf (S)}$ 
that satisfies the conditions in Theorem \ref{existence} 
be given.
Assume that for all $e\in E$ we have 
\eqref{jassumption}.
Let $T > 0$  be given.
For all $e\in E$, let initial states 
$y^e_+$,  $y^e_-$, $y^e_0$, $z^e_+$, $z^e_-$, $z^e_0
\in  C^1(0, \, L^e)$ be given
such that at all nodes $ v\in V$  the $C^1$-compatibility conditions
are satisfied.
Assume that  ${\bf (S) }$ and ${\bf (R)}$
have classical solutions on $[0, T]$
such that 
(\ref{lunendlichboundpm})
and the following
a-priori bounds hold:
%
%
%
%
%
%
There exists a numbers 
$\tMz>0$, 
$c>0$ and
$\varepsilon_0>0$
such that
\eqref{aprioribound2a},
\eqref{aprioribound1},
\eqref{aprioribound2}
hold and
\begin{equation}
\label{aprioribound20}
 |\partial_x \left( \lambda_0^e( R ) \right)| \leq \varepsilon_0
. 
\end{equation}

Assume  that there is a nonempty compact  convex set $ {\cal U} \subset {\mathbb R^2}$
that contains the function values
of the solutions 
 $(S^e_+, \, S^e_-)$ and $(R^e_+, R^e_-)$ and a number $\beta>0$ such that
\eqref{aprioribound4}
holds
and we have
\begin{equation}
\label{aprioribound40}
\max_{ \tilde R \in {\cal U}  } 
\left\{ 
\left| \partial_{R_+} \lambda_0^e( \tilde R )\right|,\,
\left|\partial_{R_-} \lambda_0^e( \tilde R )\right| \right\}
\leq \beta.
\end{equation}

Assume that there exist  numbers $\textswab{v} >0 $,  
$\overline{\textswab{v}} >0$,
such that for all $e\in E$
\begin{equation}
\label{aprioribound10}
 \textswab{v} \leq  |\lambda_0( R_+^e, \, R_-^e)|  \leq \overline{\textswab{v}}.
\end{equation}
 

Then on $[0, \, T]$ there exists a classical
solution $(\delta_+, \delta_-,\, \delta_0)$
of system ${\bf (Diff)}$.
%
\\
Let $\psi \in (0, \, \infty)$ 
be given such that  for all $e\in E$ we have 
\eqref{10012023}
(with $\kappa^e$ as defined in \eqref{kappadefinition})
and

Assume that for each  node
$v\in V$  the  number 
$\mu^v\in [0, 1]$ is 
sufficiently small in the sense that
for all $t\in [0,\, T]$
inequality \eqref{inequality2023} holds.
%
%
Assume that the perturbations are bounded above by 
$\eta_0$, in the sense that for all $t\in [0,\, T]$ we have 
\eqref{pertbound} 
and
\begin{equation}
    \label{pertbound0}
3 
\sum_{v\in V}
\sum_{e \not \in E_{in}(v,t)}
\exp( -s^e\, \psi_0\, x^e(v))\,|\lambda_0^e(S + \delta)|(t,x^e(v))  
 ( 1 - \mu^v)  \left|Z_0^v(t)\right|^2
\leq \eta_0.
\end{equation}
Then  the  solution of  ${\bf (Diff)}$
%
decays   exponentially fast 
up to the perturbation level 
in the sense that
there exist  constants
$C_0 >0$, 
$\chi_0>0$
such that for all $t\in [0,\, T]$ 
we have
\begin{equation}
\label{shirilyapunov0}
\sum_{e\in E}
 \int_0^{L^e} \left|\delta^e_0(t,\, x)  \right|^2 
\, dx
\leq C_0
\exp(-\chi_0 \,  t) \,
\sum_{e\in E}
 \int_0^{L^e} 
 \left|\delta^e_0(0,\, x)  \right|^2
 + C_0\, \left[\eta_0 + \eta_{\sigma} \right].
\end{equation}
Hence the $L^2$-norm of the error $\delta_0$ between the  state $R_0$
of the observer and the  state  $S_0$ of the original system  decays exponentially fast
up to the perturbation level.
\end{theorem}
\begin{proof}
Due to
\eqref{velocityboundinitial}
Theorem \ref{existence} implies that
\eqref{signumvor} holds and we have 
$$s^e  = 
{sign}(\lambda_0(J^e_+,\,J^e_-)(x) )
=
{sign}(\lambda_0(S^e_+,\,S^e_-)(t, \, x) ).
$$
In particular, $s^e$ is independent of $t$ and $x$.
We use the notation
$$
 W^e(S+\delta)(t,x)
=
\lambda_0(S^e_+ + \delta^e_+,\,S^e_- + \delta^e_-)(t, \, x)).
$$
For all $e\in E$ we  have
$(\delta^e_0)_t
+
W^e(S+\delta)
(\delta^e_0)_x 
= - [ W^e(S+ \delta)  - W^e(S) ] \, S^e_x.
$

With integration  by parts this yields
$\frac{d}{dt} {\cal E}^e_0(t)$
\begin{eqnarray*}
&  =   &
 \int_0^{L^e} 
2 \exp( - s^e\, \psi_0\, x)
\,
 \delta^e_0\,  {\color{black}  \partial_t\delta^e_0} \,dx
 \\
&
=
&
 -\int_0^{L^e} 
2 \exp( - s^e\, \psi_0\, x) \,
 \delta^e_0\,  W^e(S+ \delta) \, {\color{black} \partial_x \delta^e_0}
 +2 \exp( - s^e\, \psi_0\, x)
 {\color{black} \delta^e_0} \, [W(S^e + \delta^e) - W(S^e)] (S^e)_x
 \,dx
 \\
 &
= 
&
 -\int_0^{L^e} 
\exp( - s^e\, \psi_0\, x) 
\, W^e(S+ \delta)\,
 \left( |\delta^e_0|^2 \right)_x 
+2 \exp( - s^e\, \psi_0\, x)
 \delta^e_0\, 
[W(S^e + \delta^e) - W(S^e)] (S^e)_x
 \,dx
 \\
&= &  
- 
\textrm{e}^{ - s^e \, \psi_0\, x} \, W^e(S+ \delta)  |\delta^e_0(t,\, x)|^2|_{x=0}^{L^e}
 +
 \int\limits_0^{L^e} 
\left[- \psi_0 
\, \left|W^e(S+ \delta)  \right| + W^e_x(S+ \delta)\right]
e^{ - s^e\,\psi_0\, x}
 |\delta^e_0|^2 
 \\
 & & 
 +
 2 \exp( - s^e \, \psi_0\, x)
 \delta^e_0\, 
 [W(S^e + \delta^e) - W(S^e)] \,  (S^e)_x  \,dx
 .
 \end{eqnarray*}
With
\eqref{aprioribound40}
we obtain
$| W(S^e + \delta^e) - W(S^e) | \leq \beta \, \left( |\delta^e_+| + |\delta^e_-|\right)$.
Using
\eqref{aprioribound10},
\eqref{aprioribound20}
and
\eqref{aprioribound2a}
this yields the inequality
 \begin{eqnarray*}
\frac{d}{dt} {\cal E}^e_0(t)
 &
\leq
&
\left[- \psi_0 \,\textswab{v}
+
\varepsilon_0 
\right]
{\cal E}^e_0(t)
+ 2\, \sqrt{2}
 \exp( \psi\, L^e)
\beta\, \tMz \,
\int\limits_0^{L^e} 
|\delta^0_0|
[
(\delta^e_+)^2 + (\delta^e_-)^2]^{1/2}
\,
dx
\\
&
-
&
\exp( - s^e\, \psi_0\, x)\,W^e(S + \delta)(t,x)  |\delta^e_0(t,\, x)|^2|_{x=0}^{L^e}
\\
&
\leq
&
\left[- \psi_0 \, \textswab{v}
+
\varepsilon_0
\right]{\cal E}^e_0(t)
+ 2\, \sqrt{2}
 \exp( \psi\, L^e)
\beta\, \tMz 
[\int\limits_0^{L^e} 
|\delta^0_0|^2\, dx]^{ \frac{1}{2} }
[\int\limits_0^{L^e} 
(\delta^e_+)^2 + (\delta^e_-)^2
\,
dx]^{\frac{1}{2}}
\\
&
-
&
\exp( -s^e\, \psi_0\, x)\,W^e(S + \delta)(t,x)  |\delta^e_0(t,\, x)|^2|_{x=0}^{L^e}
.
 \end{eqnarray*}
Define
$\chi_0 :=
\psi_0 \, \textswab{v}
-
\varepsilon_0
$,
$
{\cal E}_0(t)
:=
\sum_{e\in E}
{\cal E}^e_0(t)
$
and
$$Q_0(t):=
\sum_{e\in E} \exp( -s^e\, \psi_0\, x)\,W^e(S + \delta)(t,x)  |\delta^e_0(t,\, x)|^2|_{x=0}^{L^e}.
$$
Then with 
$D_0= 2\, \sqrt{2}
 \exp(( \psi + \psi_0) \, L^e)
\beta\, \tMz 
$
we have
\begin{equation}
\label{schranke1}
\frac{d}{dt} {\cal E}_0(t)
\leq
- \chi_0
{\cal E}_0(t)
+ D_0
\sqrt{{\cal E}_0(t)} 
\left[
\sum_{e\in E} \int\limits_0^{L^e} 
(\delta^e_+)^2 + (\delta^e_-)^2
\,
dx
\right]^{\frac{1}{2}}
+
Q_0(t).
\end{equation}
To finish the proof, a detailed analysis of the nodal term $Q_0(t)$ is necessary.
Due to \eqref{aprioribound10}
for all $v\in V$
the sets
$e\in E_{in}(v,t)$ are independent of $t$.
Hence we have $E_{in}(v,t) = E_{in}(v,0)$.
We have
\begin{multline}\label{eq:Q0} Q_0(t)=
\sum_{v\in V}\Big[
\sum_{e\in E_{in}(v,0)}
\exp( -s^e\, \psi_0\, x^e(v))\,|W^e(S + \delta)|(t,x^e(v))  |\delta^e_0(t,\, x)|^2
\\
-
\sum_{e  \in E(v) \setminus E_{in}(v,0)}
\exp( -s^e\, \psi_0\, x^e(v))\,|W^e(S + \delta)|(t,x^e(v))  |\delta^e_0(t,\, x)|^2.\Big]
\end{multline}

Note that the coupling condition for $\delta_0$ in \eqref{Diff} can be expressed as
\begin{multline}\label{diffd0}
\delta_0^e(t,\, x^e(v))  = ( 1 - \mu^v)  Z_0^v(t)
+
\mu^v \, K^v_0(\delta_0, \, S_+ + \delta_+, \, S_-+\delta_-,\, t)
\\
+
\mu^v \, [K^v_0(S_0, \, S_+ + \delta_+, \, S_-+\delta_-,\, t)
-
 K^v_0(S_0, \, S_+ , \, S_-,\, t)]
\end{multline}
since $K^v_0$ is linear in its first component.
%
Inserting \eqref{diffd0} into \eqref{eq:Q0}  yields

\begin{equation}
\begin{split}
Q_0(t) &\geq 
\sum_{v\in V} \Bigg[
\sum_{e\in E_{in}(v, 0)}
\exp( -s^e\, \psi_0\, x^e(v))\,|W^e(S + \delta)|(t,x^e(v))  |\delta^e_0(t,\, x)|^2
\\
&-
\sum_{e  \in E(v) \setminus E_{in}(v, 0)}
\exp( -s^e\, \psi_0\, x^e(v))\,|W^e(S + \delta)|(t,x^e(v))  
\\
& \ \ \times \Big\{\Big|
 ( 1 - \mu^v)  Z_0^v(t)
+
\mu^v \, K^v_0(\delta_0, \, S_+ + \delta_+, \, S_-+\delta_-,\, t)
\\
& \ \ +
\mu^v \, [K^v_0(S_0, \, S_+ + \delta_+, \, S_-+\delta_-,\, t)
-
 K^v_0(S_0, \, S_+ , \, S_-,\, t)]
\Big|^2 \Big\}\Bigg]
.
\end{split}
\end{equation}
For $v\in V$ we define 
\begin{equation*}
\begin{split}
    X_1^v &:= 
3
\sum_{e  \in E(v) \setminus E_{in}(v, 0)}
\exp( -s^e\, \psi_0\, x^e(v))\,|W^e(S + \delta)|(t,x^e(v)) \\
& \qquad \qquad \qquad \qquad \qquad 
\times \left|
\mu^v \, K^v_0(\delta_0, \, S_+ + \delta_+, \, S_-+\delta_-,\, t)
\right|^2,
\\
X_2^v & :=
3 \sum_{e  \in E(v) \setminus E_{in}(v, 0)}
\exp( -s^e\, \psi_0\, x^e(v))\,|W^e(S + \delta)|(t,x^e(v))
\\
& \qquad\qquad \quad \times \left|
\mu^v \, [K^v_0(S_0, \, S_+ + \delta_+, \, S_-+\delta_-,\, t)
-
 K^v_0(S_0, \, S_+ , \, S_-,\, t)] 
\right|^2, 
\\
X_3^v &: =
3
\sum_{e  \in E(v) \setminus E_{in}(v, 0)}
\exp( -s^e\, \psi_0\, x^e(v))\,|W^e(S + \delta)|(t,x^e(v))  
 ( 1 - \mu^v)  \left|Z_0^v(t)\right|^2.
 \end{split}
\end{equation*}
Since
for any 
$x_1,x_2,x_3 \in {\mathbb R}$ we have
$(x_1 + x_2 + x_3)^2 \leq 3 x_1^3 + 3x_2^2 + 3 x_3^2$
this yields
$$
Q_0(t) \geq 
\sum\limits_{v\in V}
\sum\limits_{e\in E_{in}(v, 0)}
{\rm e}^{  -s^e\, \psi_0\, x^e(v) } 
\,|W^e(S + \delta)|(t,x^e(v))  |\delta^e_0(t,\, x)|^2 
-X_1^v -X_2^v  - X_3^v.
$$

The definition of
$K^v_0(\delta_0, \, S_+ + \delta_+, \, S_-+\delta_-,\, t)$
implies that similarly as for $Q_{\sigma}$
 for $\mu^v$
sufficiently small,
we have
$$
\sum_{v\in V}
\sum_{e\in E_{in}(t)}
\exp( -s^e\, \psi_0\, x^e(v))\,|W^e(S + \delta)|(t,x^e(v))  |\delta^e_0(t,\, x)|^2
-X_1^v\geq 0.
$$
Hence 
$Q_0(t)$
is bounded below by a term
that is determined by
$\sum\limits_{v\in V} X_2^v$
and $\sum\limits_{v\in V} X_3^v$.
Due to \eqref{pertbound0} 
we have
$\sum_{v\in V} X_3^v \leq \eta_0$.
In Theorem \ref{l2lyapunovthm}
it is stated that
$(\delta_+, \, \delta_-)$
decays exponentially fast
and thus can be bounded above
by an exponential decaying term 
and a term of the order
${\cal O}(\eta_\sigma)$. 
The bound (\ref{velocityboundstate})
implies that there is no pipe with zero velocity.
Hence due to 
(\ref{Kirchhoff}), 
$E_{in}(v,0)$ is not empty and (\ref{aprioribound10}) implies

$
\sum_{g \in E_{in}(v,0)} (D^g)^2 |\lambda_0^g(t, \, x^g(v))|
\geq 
\textswab{v} \,\min_{g \in E_0(v)} |D^g|^2   .
$

Hence the denominators 
that appear in $\lambda^f_R(t)$
in the definition (\ref{k0definition})  of
$K^v_0$  can be strictly bounded away from zero.
Thus we obtain
a bound for 
$\sum_{v\in V} X_2^v$ of the same order
${\cal O}(\exp( -{\chi} \, t) + \eta_\sigma)$.
%
Hence by further increasing the
number $D_0>0$ if necessary 
due to
\eqref{shirilyapunov}
and  (\ref{schranke1})
we have 
 \begin{eqnarray*}
\frac{d}{dt} {\cal E}_0(t)
& \leq
&
- \chi_0 \, {\cal E}_0(t) 
+ 
D_0 \, \sqrt{{\cal E}_0(t)} 
\sqrt{\exp( -{\chi} \, t) + \eta_\sigma}
-
Q_0(t)
 \\
 & \leq &
 - \chi_0 \, {\cal E}_0(t) 
+ 
D_0 \, \sqrt{{\cal E}_0(t)} 
\sqrt{\exp( -{\chi} \, t) + \eta_\sigma}
+
\martin{
D_0
(
 \eta_\sigma
+
\exp( - {\chi}\, t))
+
 \eta_0.
 }
 \end{eqnarray*}
Therefore  \eqref{gwl2vor} holds with $\eta= \eta_\sigma$.
Then Lemma \ref{gwl2} implies the assertion.


\end{proof}









\section{Conclusion}
We have defined an observer system for the gas flow through a pipeline network 
that is governed by a quasilinear model.
The system allows for hydrogen blending in the natural gas flow.
As input into the observer system, 
measurement data that is  obtained at 
the nodes 
of  the network is used.
We have allowed for a certain measurement error
with the assumption that is is 
smoothed to have $C^1$ regularity.

We have shown that under suitable
regularity conditions for the solution
the observation error decays exponentially fast
up to the level of the measurement error.
In the proofs of our results appropriately chosen
Lyapunov functions 
with exponential weights play an essential role.
Our result requires smallness assumptions on the initial
error.
To be precise, the initial data of both the original and the observer system need to be close to a
given 
steady reference state.
It  remains  a challenge 
for future studies
to find results that allow for large initial
errors.

\textbf{Acknowledgements:}
This work was
supported by DFG in the 
 Collaborative Research Centre
CRC/Transregio 154,
Mathematical Modelling, Simulation and Optimization Using the Example of Gas Networks,
Project C03 and  Project  C05.
The authors thank the Bundesministerium für Bildung und Forschung (BMBF)
for support under DAAD grant 57654073 'Uncertain data in control of PDE systems'.

\bibliographystyle{siamplain}

\end{document}